\documentstyle[11pt,amssymb,amstex]{amsart}

\numberwithin{equation}{section}

\newcommand{\bea}{\begin{eqnarray}}
\newcommand{\eea}{\end{eqnarray}}
\newcommand{\ba}{\begin{array}}
\newcommand{\ea}{\end{array}}
\newcommand{\edc}{\end{document}}
\newcommand{\bc}{\begin{center}}
\newcommand{\ec}{\end{center}}
\newcommand{\be}{\begin{equation}}
\newcommand{\ee}{\end{equation}}

\def\bb{{\mathbb B}}
\def\bc{{\mathbb C}}

\def\bn{{\mathbb N}}
\def\bp{{\mathbb P}}

\def\br{{\mathbb R}}

\def\a{\alpha}
\def\b{\beta}
  
\def\d{\delta}  

\def\e{\varepsilon}
\def\l{\lambda} 
\def\k{\kappa}
\def\m{\mu}

\def\s{\sigma} 
\def\t{\tau}

\def\th{\theta}

\def\ab{{\mathbf{a}}}
\def\id{{\bf 1}\!\!{\rm I}}
\newtheorem{thm}{Theorem}[section]
\newtheorem{lem}[thm]{Lemma}
\newtheorem{cor}[thm]{Corollary}

\begin{document}
\small

\title[Local ergodic theorem]
{On noncommutative weighted local ergodic theorems on
$L^p$-spaces}

\author{Farrukh Mukhamedov}
\address{Farrukh Mukhamedov\\
 Department of Computational \& Theoretical Sciences\\
Faculty of Sciences, International Islamic University Malaysia\\
P.O. Box, 141, 25710, Kuantan\\
Pahang, Malaysia} \email{{\tt far75m@@yandex.ru} {\tt
farruh@@fis.ua.pt}}
\author{Abdusalom Karimov}
\address{Abdusalom Karimov\\
Department of Mathematics\\
Tashkent Institute of Textile and Light Industry\\
Tashkent, 700100, Uzbekistan} \email{\tt karimov57@@rambler.ru}

\begin{abstract}

In the present paper we consider a von Neumann algebra $M$ with a
faithful normal semi-finite trace $\t$, and $\{ \alpha_ t\} $  a
strongly continuous extension to $L^p(M,\t)$ of a semigroup of
absolute contractions on $L^1 (M,\tau )$. By means of a
non-commutative Banach Principle we prove for a Besicovitch
function $b$ and $x\in L^p(M,\t)$, the averages
\begin{equation*}
\frac{1}{T}\int_0^Tb(t)\a_t(x)dt
\end{equation*}
converge bilateral almost uniform in $L^p(M,\t)$ as $T\to 0$.

\end{abstract}

\maketitle

\section{Introduction}

It is known (see for example \cite{k}) that in the classical
ergodic theory one of the powerful tools in dealing with the
almost everywhere convergence of ergodic avarages is the well-
known Banach Principle, which can be formulated as follows:

\begin{thm} Let $(S, F,m)$ be a measurable space
with a $\s$-finite measure and let $X$ be a Banach space. Let
$\{a_n\}$ be a sequence of continuous linear maps of $X$ into the
space of measurable functions on $S$. Assume that
$\sup_n\{|a_n(x)(s)|\} < \infty$ for each $x\in X$ and almost all
$s\in S$. If the sequence $a_n(x)$ converges almost everywhere for
$x$ in a dense subset of $X$, then this sequence converges for
each $x\in X$.
\end{thm}

This principle is often applied in proofs concerning the almost
everywhere convergence of weighted averages, moving averages, etc.

In a non-commutative setting the almost everywhere convergence of
sequences of operators were applied to study of the individual
ergodic theorems in von Neumann algebras by many authors
\cite{g},\cite{gg},\cite{La},\cite{p},\cite{y} (see \cite{ja1} for
review). But in these investigations those ergodic theorems were
obtained without using an analog of the Banach Principle. In
\cite{gl} firstly a non-commutative analog of such principle was
proved for quasi-uniform convergence. Using that result in
\cite{lm} a uniform sequence weighted ergodic theorem was proved
in the space of integrable operators affiliated with a von Neumann
algebra. Recently, in \cite{CLS} for the Banach Principle for
bilateral uniform convergence has been adopted, and  by means of
it the Besicovitch weighted ergodic theorem has been proved.

In the present paper we are going to prove local and weighted
local ergodic theorems on non-commutative $L^p$-spaces by means of
the Banach principle. Note that such kind of theorems in
commutative settings were studied by many authors (see for example
\cite{af},\cite{hs},\cite{k}). In a non-commutative setting we
mention works \cite{c}, \cite{cd}, \cite{ja1}, \cite{jx2},
\cite{w}.

Let us end this section with description of the organization of
the paper.  In Section 2, we recall some preliminary results and
formulate the Banach Principle. There, to prove local ergodic
theorem, we adopt the principle in a more convenient formulation.
In the next Section 3 we prove the local ergodic theorem for
semigroups of absolute contractions of $L^p$-spaces. Note that
this section reviews the results of \cite{c},\cite{jx2}.  Using
the result of Sec. 3, in  the last Section 4, we establish a
weighted local ergodic theorem by means of the Banach principle.

\section{Preliminaries}

Let $M$ be a semifinite von Neumann algebra acting on a Hilbert
space $H$, let $\tau$ be a  faithful normal semifinite trace on
$M$,  let $P(M)$ be the complete lattice of all projections in
$M$. A densely-defined closed operator $x$ in $H$ is said to be
{\it affiliated} with $M$ if $y^{\prime}x\subset xy^{\prime}$ for
every $y^{\prime}\in M^{\prime}$, where $M^{\prime}$ is the
commutant of the algebra $M$.  An operator $x$, affiliated with
$M$, is said to be {\it $\tau$-measurable} if for each $\e >0$
there exists $e\in P(M)$ with $\tau (e^{\perp})\leq \e$ such that
$eH\subset D_{x}$, where $e^{\perp}= \id - e $, $\id$ is the unit
of $M$, $D_{x}$ is the domain of definition of $x$. Let $S(M)$ be
the set of all $\tau-$measurable operators affiliated with $M$.
Let $\| \cdot \|$ stand for the uniform norm in $M$. The {\it
measure topology} in $S(M)$ is given by the system $$V(\e,
\d)=\{x\in S(M): \|xe\| \leq \d \text { for some } e\in P(M) \text
{ with } \tau (e^{\perp})\leq \e \}, $$ $\e >0$, $\d >0$, of
neighborhoods of zero. Accordingly, a sequence $\{x_n\} \subset
S(M)$ converges {\it in measure} to $x\in S(M)$, $x_n\to x$ (m),
if, given $\e>0, \d>0$, there is a number $N=N(\e,\d)$ such that
for any $n\geq N$ there exists a projection $e_n \in P(M)$
satisfying the conditions $\tau (e_n^\perp)<\e$ and $\| (x_n-x)e_n
\|<\d$.

\begin{thm}\cite{Ne} Equipped with the
measure topology, $S(M)$ is a complete topological *-algebra.
\end{thm}

For a positive self-adjoint operator $x=\int_{0}^{\infty}\lambda
de_{\lambda}$ affiliated with $M$ one can define
$$\tau(x)=\sup_{n}\tau \left (\int_{0}^{n}\lambda de_{\lambda}\right )
=\int_{0}^{\infty} d\tau(e_{\lambda}).$$
 If $0<p\leq \infty$, then
$$L^{p}=L^{p}(M, \tau)=\cases \{x\in S(M): \|x\|_{p}=\tau
(\arrowvert x\arrowvert^{p})^{1/p}<\infty \},&\text{for $p\ne
\infty$}\\ (M, \|\cdot \| ),&\text{for $p=\infty$}.\endcases$$
Here, $\arrowvert x\arrowvert$ is the {\it absolute value} of $x$,
i.e. the square root of $x^{*}x$. By $L^p_+$ ( resp. $L^p_{sa}$)
we denote the set of positive (resp. self-adjoint) elements of
$L^p$.  We refer a reader to \cite{px} for more information about
noncommutative integration and to \cite{sa,t} for general
terminology of von Neumman algebras.

There are several different types of convergences in $S(M)$ each
of which, in the commutative case with finite measure, reduces to
the  almost everywhere convergence (see for example \cite{pa}). In
the paper we deal with so called the {\it bilateral almost uniform
} (b.a.u.) convergence in $S(M)$ for which $x_{n}\to  x$ means
that for every $\e>0$ there exists $e\in P(M)$ with
$\tau(e^{\perp})\leq \e$ such that $\|e(x_n- x)e\|\to 0$. It is
clear that b.a.u. implies convergence in measure.  Now recall well
known fact concerning b.a.u. convergence (see \cite{Ne},\cite{s}).

\begin{lem}\label{ad} Let $M$ be as above. If two sequences $x_n$
and $y_n$ converge b.a.u., then $x_n+y_n$ converges b.a.u.
\end{lem}

In \cite{CLS} the following results has been proved.

\begin{thm}\label{comp} Algebra $S(M)$ is complete
with respect to the b.a.u. convergence.
\end{thm}

\begin{lem}\label{lp}
Let $0\leq p<\infty$, and let $\{ x_n\} \subset L^p$ be such that
$\liminf_{n} \| x_n\|_p=s<\infty$. If $x_n \to \ x$ b.a.u., then
$\ x\in L^p$ and $\| x\|_p\leq s$.
\end{lem}

Recall a non-commutative the Banach Principle (see \cite{CLS}).
Let $(X, \ \| \cdot \|, \ \geq )$ be an ordered real Banach space
with the closed convex cone $X_+$, $X=X_+-X_+$. A subset $X_0
\subset X_+$ is said to be {\it minorantly dense} in $X_+$ if for
every $x\in X_+$ there is a sequence $\{ x_n \}$ in $X_0$ such
that $x_n \leq x$ for each $n$, and $\| x-x_n \| \to 0$ as $n\to
\infty$. A linear map $a:X\to S(M)$ is called {\it positive} if
$a(x)\geq 0$ whenever $x\in X_+$.

\begin{thm}\label{bp} Let $X$ be an ordered real Banach space
with the closed convex cone $X_+$. Let $a_{n}:X \to S(M)$ be a
sequence of positive continuous (in the measure topology)  linear
maps satisfying the conditions
\begin{itemize}
    \item[(i)] For every $x \in X_+$ and $\e>0$ there is $b\in M$, $0\neq
b\leq I$, such that $\tau (I-b)<\e$, and
$$\sup_{n} \| ba_{n}(x)b \|<\infty.$$
    If, for every $x$ from a minorantly dense subset $X_{0}\subset
X_+$, \item[(ii)] $a_{m}(x)-a_{n}(x) \to 0 \ \ $ b.a.u., $m,n\to
\infty$,
\end{itemize}
then (ii) holds on all of $X$.
\end{thm}

{\bf Remark.} According to Theorem \ref{comp} that  the
fundamental sequences in Theorem \ref{bp} indeed have their limits
belonging to $S(M)$.

As it has been pointed out that the Banach Principle is one of the
basic tools to prove ergodic theorems. But the above formulated
Principle is too complicated to apply, since it requires
minorantly density of $X_0$, which makes difficult to check the
condition (ii). Basically, to obtain some ergodic theorems we
really need the following theorem, which is an analog of the
Banach Principle.

\begin{thm}\label{bp2} Let $X$ be a Banach space and let $\ab_{n}:X\to S(M)$ be a
sequence of linear maps satisfying the conditions
\begin{itemize}
    \item[(i)] For every $x \in X$ and $\e>0$ there is $p\in P(M)$, with $\tau (p^{\perp})<C(\e^{-1}\|x\|_X)^\a$,
    such that  $\| p ~ \ab_{n}(x)p\|<\e$ for all
    $n\in\bn$, here $C$ and $\a$ are some positive constants.\\
   If, for every $x$ from a dense subset $X_{0}\subset
X$, \item[(ii)] $\ab_{m}(x)-\ab_{n}(x) \to 0 \ \ $ b.a.u., $m,n\to
\infty$,
\end{itemize}
then (ii) holds on all of $X$.
\end{thm}

\begin{pf} Let $x\in X$. Due to density of $X_0$ in
$X$, for given $\e>0$ there is a sequence $\{x_n\}\subset X_0$
such that $\|x_n-x\|_X<(\e/2^{n+1})^{2/\a}$ for every $n\in\bn$.
Then from (i) for every $n\in\bn$ there is a projection $p_n\in
P(M)$ with $\t(p_n^{\perp})<C\e/2^{n+1}$ such that
\begin{equation}\label{1m}
\|p_n(\ab_m(x_n-x))p_n\|<\e/2^{n+1} \ \qquad \forall m\in\bn.
\end{equation}

Putting $p=\bigwedge\limits_n p_n$, we have $\t(p^{\perp})<C\e/2$
and
\begin{equation*}
\|p(\ab_m(x_n-x))p\|\to 0, \ \ n\to\infty \ \ \textrm{uniformly in
} \ m.
\end{equation*}

Therefore, from the last relation for given $\e>0$ one finds
$n_0\in\bn$ such that
\begin{equation}\label{2m}
\|p(\ab_m(x_{n_0}-x))p\|\leq \frac{\e}{3}
\end{equation}
for all $m\in\bn$. Since $x_{n_0}\in X_0$ by condition (ii) there
is a projection $q\in P(M)$ with $\t(q^{\perp})<\e/2$ and
$N_0\in\bn$ such that
\begin{equation}\label{3m}
\|q(\ab_{m}(x_{n_0})-\ab_{n}(x_{n_0}))q\|\leq \frac{\e}{3}
\end{equation}
for all $m,n\geq N_0$. Letting $f=p\wedge q$ one gets
$\t(f^{\perp})<\e(C+1)/2$ and \eqref{2m},\eqref{3m} imply
\begin{eqnarray*}
\|f(\ab_m(x)-\ab_n(x))f\|&\leq&
\|p(\ab_m(x_{n_0}-x))p\|+\|p(\ab_n(x_{n_0}-x)p\|\\
&&+\|q(\ab_{m}(x_{n_0})-\ab_{n}(x_{n_0}))q\|\leq \e
\end{eqnarray*}
This proves the assertion.
\end{pf}

{\bf Remark.} We should note that in the proved Theorem a Banach
space $X$ need not be ordered. Hence a condition of minorantly
density of $X_0$ and positivity of $\ab_n$ are extra restrictions,
which were important in Theorem \ref{bp}. But the condition (i) in
Theorem \ref{bp2} is strong than one in Theorem \ref{bp}. For
example, it implies that each mapping $\ab_m$ ($m\in\bn$) is
continuous with respect to b.a.u. convergence, which can be seen
from \eqref{1m}.

Recall a positive linear map $\a :L^1(M,\tau)\to L^1(M,\tau)$ will
be called an {\it absolute contraction} if $\a(x)\leq \id$ and
$\tau(\a(x))\leq \tau(x)$ for every $x\in M\cap L^{1}$ with $0\leq
x\leq \id$. If $\a$ is a positive contraction in $L^1$, then, as
it is shown in \cite{y}, $\|\a(x)\|_{p}\leq \|x\|_{p}$ holds for
each $x=x^*\in M\cap L^p$ and all $1\leq p\leq \infty$. Besides,
there exist unique continuous extensions $\a :L^{p}\to L^{p}$ for
all $1\leq p<\infty$ and a unique ultra-weakly continuous
extension $\a : M\to M$ (see \cite{jx2},\cite{y}) . This implies
that, for every $x\in L^{p}$ and any positive integer $k$, one has
$\|\a^{k}(x)\|_{p}\leq 2\|x\|_{p}$.

Let $\{ \alpha_ t\}_{t\geq 0} $ be semigroup of absolute
contraction on $L^1 $. This means that each $\alpha_t$ is an
absolute contraction on $L^1$, $\a_0=Id$ and $\alpha_t \alpha_s
=\alpha_{t + s}$ for all $t,s \geq 0$. By the same symbol $\a_t$
we will denote its extension to $L^p$ $(1\leq p<\infty$). In the
sequel we assume that the semigroup $\{a_t\}$ is strongly
continuous in $L^p$, for fixed $p$, i.e. $\mathop {\lim
}\limits_{t \to s} \left\| {\alpha _t f - \alpha _s f} \right\|_p
= 0$ for all $s \geq 0$ and $f\in L^p$.

For any $T>0$ put
$$
\beta _T(x) = {1 \over T} \int\limits_0^T \alpha_t(x)dt \ \ \mbox{
for } \ \  x \in L^p (M,\tau ).
$$
It is clear that $\beta _T$ is  positive linear map, and maps $
L^p $ into itself. The following maximal theorem was proved in
\cite{y},\cite{jx2}.

\begin{thm}\label{max} Let $x\in L^p_{sa}$ then for any
$\varepsilon> 0$, there exists projection $e\in P(M)$ such that
$\tau(e^{\perp})< C( \varepsilon^{ - 1} \|x \|_p)^p$ and
$$
\left\| {e\beta _T (A)e} \right\|\leq \varepsilon  \ \ \textrm{
for all } \ \ T > 0.
$$
\end{thm}

\section{Local ergodic theorem}

This section is devoted to the local ergodic theorem, which can be
formulated as follows

\begin{thm}\label{loc} Let  $\{ \alpha_ t\}_{t\geq 0} $ be a strongly
continuous extension to $L^p(M,\t)$ of a semigroup of absolute
contractions on $L^1 (M,\tau )$. Then for every $x\in L^p(M,\t)$
the averages $\b_T(x)$ converge b.a.u. in $L^p(M,\t)$ as $T\to 0$.
\end{thm}

Such kind of theorems were proved in \cite{c}, \cite{cd},
\cite{ja1}, \cite{jx2}, \cite{w}. Here we are going to provide a
different proof based on the Banach Principle.

To prove the theorem, we need some auxiliary facts.

\begin{lem}\label{in} Let $x\in L^p_+$, then
\begin{equation}\label{int}
-\frac{1}{b}\int_0^a\a_s(x)ds\leq
\b_a(\b_b(x))-\b_b(x)\leq\frac{1}{b}\int_b^{b+a}\a_s(x)ds
\end{equation}
for every $a,b\in\br_+$.
\end{lem}

\begin{pf} Denote
$$
y=\int_0^b\a_s(x)ds.
$$
Then for a positive number $0<h<a$ we have
\begin{eqnarray*}
\a_h(y)-y&=&\int_h^{b+h}\a_s(x)ds-\int_0^b\a_s(x)ds\\
&=&\int_b^{b+h}\a_s(x)ds-\int_0^a\a_s(x)ds\\
&\leq&\int_b^{b+h}\a_s(x)ds
\end{eqnarray*}
here we have used that $\int_0^a\a_s(x)ds\geq 0$. Whence
\begin{eqnarray}\label{1}
\b_a(y)-y&=&\frac{1}{a}\int_0^a(\a_h(y)-y)dh\nonumber\\
&\leq&\frac{1}{a}\int_0^a\bigg(\int_b^{b+h}\a_s(x)ds\bigg)dh\nonumber\\
&\leq&\frac{1}{a}\int_0^a\bigg(\int_b^{b+a}\a_s(x)ds\bigg)dh\nonumber\\
&=&\int_b^{b+a}\a_s(x)ds.
\end{eqnarray}
The last inequality \eqref{1} implies
\begin{equation}\label{2}
\b_a(\b_b(x))-\b_b(x)\leq\frac{1}{b}\int_b^{b+a}\a_s(x)ds.
\end{equation}
On the other hand, we have
\begin{eqnarray*}
\a_h(y)&=&\int_h^{b+h}\a_s(x)ds\\
&\geq& \int_h^b\a_s(x)ds\\
&\geq&\int_a^b\a_s(x)ds  \ \ \textrm{for} \ \ 0<h<a.
\end{eqnarray*}
Therefore,
$$
\int_a^b\a_s(x)ds\leq\frac{1}{a}\int_0^a\a_h(y)dh
$$
which yields
$$
-\frac{1}{b}\int_0^a\a_s(x)ds\leq \b_a(\b_b(x))-\b_b(x).
$$
This and \eqref{2} complete the proof.
\end{pf}

Denote
\begin{equation}\label{x0}
X_0=span\{\b_T(x): \ x\in L^p_+, T>0\}.
\end{equation}

\begin{lem}\label{x01} The space $X_0$ is dense in $L^p$.
\end{lem}

\begin{pf} Take $x\in L^1$, and show there is a sequence $\{x_k\}$ in
$X_0$  which converges to $x$ in norm of $L^1$. Define a sequence
$\{x_k\}$ by
\begin{equation}\label{xk}
x_k=k\int_0^{1/k}\a_s(x)ds.
\end{equation}
Since any $x \in L^1$ can be represented by $x
=\sum\limits_{j=0}^3 i^kx_j$, where $x_j\in L^p _{+}$ (
$j=0,1,2,3$), therefore $x_k$ is a linear combination of
$\b_{1/k}(x_j)$, which implies that $\{x_k\}\subset X_0$.  The
strong continuity of $\a_s$ implies that for arbitrary $\e>0$
there is $\d>0$ such that for every $s$ with $|s|<\d$ the
inequality holds $\|\a_s(x)-x\|_p<\e$. Pick $k_0\in\bn$ such that
$k_0<\d$, then
\begin{equation*}
\|x_n-x\|_p\leq n\int_{0}^{1/n}\|\a_s(x)-x\|_pds<\e \ \ \ \forall
n\geq k_0
\end{equation*}
which completes the proof.
\end{pf}

\begin{lem}\label{bb} Let $x\in L^p_+$, then
\begin{equation}\label{bb1}
\lim_{a\to 0}\b_a(\b_b(x))=\b_b(x) \ \ b.a.u.
\end{equation}
for every $b>0$.
\end{lem}

\begin{pf} First denote
\begin{eqnarray}\label{bb2}
h(a)=\frac{1}{b}\int_0^a\a_s(x)ds, \qquad
g(a)=\frac{1}{b}\int_b^{b+a}\a_s(x)ds,
\end{eqnarray}
it is obvious that $h(a)\geq 0$, $g(a)\geq 0$ for all $a>0$. Now
due to the strong continuity of $\a_s$ we infer
\begin{eqnarray*}
\lim_{a\to 0}\|h(a)\|_p=0, \qquad \lim_{a\to 0}\|g(a)\|_p=0.
\end{eqnarray*}
From this we conclude that for any $\e>0$ there is a sequence
$\{a_k\}\subset\br_+$ such that $\t(h^p(a_k))<\e^2/2^{2k}$ for all
$k\in\bn$.

Let
$$
h^p(a_k)=\int_0^\infty\l de^{(k)}_\l
$$
be the spectral resolution of $h^p(a_k)$. Put
$p_k=e^{(k)}_{\e/2^{k+1}}$, then $\t(p_k^{\perp})\leq \e/2^{k+1}$.
From
$$
p_kh(a_k)p_k=\int_0^{\e/2^{k+1}}\l^{1/p} de^{(k)}_\l
$$
one sees that $p_kh(a_k)p_k\in M$, and with the inequality
$h(a)\leq h(c)$ for $0<a<c$ for sufficiently small $a$ we have
\begin{eqnarray}\label{bb2}
\|p_kh(a)p_k\|\leq\|p_kh(a_k)p_k\|\leq \frac{\e}{2^{k+1}}.
\end{eqnarray}
Letting $p=\bigwedge p_k$, one finds $\t(p^{\perp})<\e/2$. It
follows from \eqref{bb2} that
\begin{eqnarray}\label{bb3}
\|ph(a)p\|\leq\|p_kh(a)p_k\|\leq \frac{\e}{2^{k+1}} \qquad
\textrm{for all} \ \ k\in\bn.
\end{eqnarray}

By the same argument one finds $q\in P(M)$ with
$\t(q^{\perp})<\e/2$ such that
\begin{eqnarray}\label{bb4}
\|q g(a)q\|\to 0 \qquad \textrm{as} \ \ a\to 0.
\end{eqnarray}

Put $e=p\wedge q$, then $\t(e^{\perp})<\e$. Now Lemma \ref{in}
implies that
$$
eh(a)e\leq e(\b_a(\b_b(x))-\b_b(x))e\leq eg(a)e
$$
whence from \eqref{bb3}-\eqref{bb4} one gets
$$
\|e(\b_a(\b_b(x))-\b_b(x))e\|\leq\max\{\|eh(a)e\|,\|eg(a)e\|\}\to
0 \qquad \textrm{as} \ \ a\to 0.
$$
This competes the proof. \end{pf}

The proved lemma and Lemma \ref{ad} yields the following
\begin{cor}\label{bbc} For any $x\in X_0$, we have
\begin{equation*}
\lim_{a\to 0}\b_a(x)=x \ \ b.a.u.
\end{equation*}
\end{cor}

Now we are ready to prove the formulated Theorem \ref{loc}.

\begin{pf} Take $X=L^p$ in Theorem \ref{bp2}. Then due to Theorem \ref{max} the condition (i) of
Theorem \ref{bp2} is satisfied. Now take an arbitrary sequence of
positive numbers $\{a_n\}$ such that $a_n\to 0$. Then according to
Lemma \ref{bbc} one sees that $\b_{a_n}(x)$ converges b.a.u. for
every $x\in X_0$. From Lemma \ref{x01} we already knew that $X_0$
is dense in $L^p$. Hence, all the conditions of Theorem \ref{bp2}
are satisfied, which implies the assertion of the theorem.
\end{pf}

{\bf Remark.} Note that similar results were proved in \cite{c}
and \cite{jx2}, respectively in $L^1$ and $L^p$ spaces. But our
approach uses the Banach principle.

\section{A weighted local ergodic theorem}

In this section by means of Theorem \ref{loc} and the Banach
principle we are going to prove a weight local ergodic theorem.

Recall that a function $P:\br_+\to\bc$ is called {\it
trigonometric polynomial} if it has the following form
\begin{equation}\label{pp}
P(t)=\sum_{j=1}^n\k_je^{2\pi i\th_jt}, \qquad t\in\br_+
\end{equation}
for some $\{\k_j\}\subset\bc$, and $\{\th_j\}\subset\br$.
 By
$\bp(\br_+)$ we denote the set of all trigonometric polynomials
defined on $\br_+$. We say that a measurable function
$b:\br_+\to\bc$ is a {\it Besicovitch function} if
\begin{itemize}
    \item[(i)] $b\in L^{\infty}(\br_+)$;
    \item[(ii)] Given any $\e>0$ there is $P\in\bp(\br_+)$ such
    that
    \begin{equation}\label{bf}
    \limsup_{T\to 0}\frac{1}{T}\int_0^T|b(t)-P(t)|dt<\e.
    \end{equation}
\end{itemize}

{\bf Remark.} A similar notion of Besicovitch weights was
introduced, for example, in \cite{dj}.

The next simple lemma will be used in the proof of main result
which was proved in \cite{lm}.

\begin{lem}\label{fun} If a  sequence $\{\tilde {a}_n\}$
in $M$ is such that for every $\e>0$ there exist a b.a.u.
convergent sequence $\{a_n\}\subset M$ and a positive integer
$n_0$ satisfying $\|\tilde {a}_n-a_n\|<\e$ for all $n\geq  n_0$,
then $\{\tilde{a}_n\}$ also converges b.a.u.
\end{lem}

The main result of this section is the following

\begin{thm}\label{wloc}
Let $M$ be a von Neumann algebra with a faithful normal
semi-finite trace $\t$, and $\{ \alpha_ t\}_{t\geq 0} $ be a
strongly continuous extension to $L^p(M,\t)$ of a semigroup of
absolute contractions on $L^1 (M,\tau )$. If $b$ is a Besicovitch
function and $x\in L^p(M,\t)$, then the averages
\begin{equation}\label{wa}
\tilde\b_T(x)=\frac{1}{T}\int_0^Tb(t)\a_t(x)dt
\end{equation}
converge b.a.u. in $L^p(M,\t)$.
\end{thm}

\begin{pf} Let $\bb$ be the unit circle in
$\bc$, i.e. $\bb=\{z\in\bc: \ |z|=1\}$. By $\m$ we denote the
normalized Lebesgue measure on $\bb$.  Let $\tilde M=M\otimes
L^\infty(\bb,\m)$ with $\tilde\t=\t\otimes\m$. Let $\tilde
L^q=L^q(\tilde M,\tilde\t)$ where $q\geq 1$.

Let us fix $\l\in\bb$ and define a map $\tilde\a^{(\l)}_t$ on
$\tilde L^1$ by
\begin{equation}\label{at}
(\tilde\a^{(\l)}_t(f))(z)=\a_t(f(\l^tz)), \qquad f\in \tilde L^1,
\ z\in\bb,\ t>0.
\end{equation}

One can see that for $f\in\tilde L^1_+$
\begin{eqnarray*}
(\tilde\a^{(\l)}_t(f))&=&\int_\bb\t(\a_t(f(\l^tz)))d\m(z)\\
&\leq&\int_\bb\t(f(\l^tz))d\m(z)=\tilde\t(f)
\end{eqnarray*}
and $\tilde\a^{(\l)}_t(\id)\leq\id$. These mean that
$\{\tilde\a^{(\l)}_t\}$ is a semigroup of absolute contractions of
$\tilde L^1$. By the same symbol denote its extension to $\tilde
L^p$. Strong continuity of $\a_t$ on $L^p$ implies that
$\tilde\a_t$ is so on $\tilde L^p$. Therefore, according to
Theorem \ref{loc} for every $f\in\tilde L^p$ the averages
$$
\frac{1}{T}\int_0^T\tilde\a^{(\l)}_t(f)dt
$$
converge b.a.u. in $\tilde L^p$ as $T\to 0$. By Lemma 4.1
\cite{CLS} we infer that the averages
$$
\frac{1}{T}\int_0^T(\tilde\a^{(\l)}_t(f))(z)dt=
\frac{1}{T}\int_0^T\a_t(f(\l^tz)dt
$$
converge b.a.u. in $L^p(M,\t)$ for almost all $z\in \bb$. Applying
this to the function $f(z)=zx$, here $x\in L^p_{+}(M,\t)\cap M$ we
obtain b.a.u. convergence of
$$
z\frac{1}{T}\int_0^T\l^t\a_t(x)dt \qquad \textrm{for almost all} \
\ z\in\bb.
$$
This implies that the averages
\begin{equation}\label{w1}
\frac{1}{T}\int_0^T\l^t\a_t(x)dt \ \ \textrm{converge b.a.u. as
$T\to 0$ for every $\l\in\bb$}.
\end{equation}

Now pick an arbitrary $\e>0$. Since $b$ is a Besicovitch function,
then there exists $P_\e\in\bp(\br_+)$ such that
$P_\e(t)=\sum_{j=1}^n\k_j\l_j^t$ and \eqref{bf} is satisfied,
where $\{k_j\}_{j=1}^n\subset\bc$, $\{\l_j\}\subset\bb$.
Consequently, from \eqref{w1} and Lemma \ref{ad} we obtain that
\begin{equation}\label{w2}
\frac{1}{T}\int_0^TP_\e(t)\a_t(x)dt
\end{equation}
converge b.a.u. as $T\to 0$.

On the other hand, from \eqref{bf} one gets
\begin{equation}\label{w3}
\bigg\|\frac{1}{T}\int_0^Tb(t)\a_t(x)dt-\frac{1}{T}\int_0^TP_\e(t)\a_t(x)dt\bigg\|=
2\bigg(\frac{1}{T}\int_0^T|P_\e(t)-b(t)|dt\bigg)\|x\|<2\e\|x\|
\end{equation}

Now Lemma \ref{fun} implies that the average \eqref{wa} converges
b.a.u. in $L^p_{+}\cap M$ as $T\to 0$. This means that b.a.u.
convergence of
\begin{equation}\label{wa1}
(\tilde\b_T(x))^*=\frac{1}{T}\int_0^T\overline{b(t)}\a_t(x)dt.
\end{equation}
The last relation with \eqref{wa} yields that both
$$
\tilde\b^{(r)}_T(x)=\frac{1}{T}\int_0^T\Re(b(t))\a_t(x)dt \qquad
\textrm{and} \qquad \tilde\b^{(i)}_T(x)=
\frac{1}{T}\int_0^T\Im(b(t))\a_t(x)dt
$$
averages converge b.a.u. too.

Put
$$
\tilde\b^{(R)}_T(x)=\tilde\b^{(r)}_T(x)+\b_T(x), \ \
\tilde\b^{(I)}_T(x)=\tilde\b^{(i)}_T(x)+\b_T(x),
$$ here as before
$$
\b_T(x)=\frac{1}{T}\int_0^T\a_t(x)dt.
$$

Now according to Theorem \ref{max} given $\e>0$ there exists a
projection $e\in P(M)$ with $\t(e^{\perp})<C(\e^{-1}\|x\|_p)^p$
such that
$$
\sup_T\|e\b_T(x)e\|<\e
$$

Note that, since $b$ from $L^\infty(\br_+)$ without loss of
generality we may assume that $|b(t)|\leq 1$ for almost every
$t\in\br_+$. Therefore, one finds $0\leq \Re (b) +1\leq 2$ which
implies that
\begin{eqnarray*}
e\tilde\b^{(R)}_T(x)e\leq 2e\b_T(x)e
\end{eqnarray*}
for every $T\in\br_+$. This immediately  yields
\begin{eqnarray*}
\|e\tilde\b^{(R)}_T(x)e\|&\leq 2\e
\end{eqnarray*}

Since $\tilde\b_T^{(R)}:X=L^p_{+}\to S(M)$ is a positive linear
continuous maps, and the set $X_0:=L^p_+\cap M$ is dense in
$X=L^p_{+}$, by Theorem \ref{bp2} we obtain the b.a.u. convergence
of $\tilde\b_T^{(R)}(x)$ for all $x\in L^p_{+}$. Remembering that
the averages $\b_T(x)$ also converge b.a.u. one gets the
convergence of $\tilde\b_T^{(r)}(x)$, $x\in L^p_{+}$. Analogously,
$\tilde\b_T^{(i)}(x)$ converges b.a.u. for all $x\in L^p_{+}$.
Therefore, by Lemma \ref{ad} the averages
$$
\tilde\b_T(x)=\tilde\b_T^{(r)}(x)+i\tilde\b_T^{(i)}(x)
$$ converge b.a.u. for every $x\in L^p_+$,
hence for every $x\in L^p$.

It remains to show that the limits of these averages belong to
$L^p$.  Taking into account that $\|\ a_t(x)\|_p\leq 2\| x\|_p$
for each $t$, we get $\|\tilde\b_T(x)\|_p\leq 2\| x\|_p$ for all
$x\in L^p$. This finishes proof due to Lemma \ref{lp}.
\end{pf}

{\bf Remark.} In the proof we could use Theorem \ref{bp} instead
of Theorem \ref{bp2}, since in that case we may take $X=L^p_{sa}$
$X_0=L^p_+\cap M$. Indeed, $X$ is an ordered Banach space with
closed cone $X_+=L^p_+$, and $X_0$ is a minorantly dense subset of
$L^p_+$ (see \cite{CLS}).

\section*{acknowledgments} The authors would like to
thank  Prof. V.I. Chilin from National University of Uzbekistan,
for valuable advice on the subject.

\end{document}